\documentclass{amsart}
\usepackage{amssymb,amsmath,amsthm}
\usepackage[shortlabels]{enumitem}
\usepackage{comment}
\usepackage{url}
\usepackage{xcolor}
\newtheorem{theorem}{Theorem}[section]

\newtheorem*{theorem*}{Theorem}

\theoremstyle{definition}

\newtheorem{conjecture}[theorem]{Conjecture}

\title{On a selfless path}

\author[S. Kunnawalkam Elayavalli]{Srivatsav Kunnawalkam Elayavalli}
\address{\parbox{\linewidth}{Department of Mathematics, University of California, Berkeley}}
\email{srivatsav.k.e@berkeley.edu}
\urladdr{https://sites.google.com/view/srivatsavke}

\begin{document}
\begin{abstract}
I give a brief account of my recent research direction concerning free independence in $C^*$-algebras i.e, \emph{selflessness}.
\end{abstract}

	\maketitle

    \subsection*{A preliminary remark}

The last few months have witnessed the communication of an unreasonable amount of important mathematical results, especially in the theory of operator algebras (see for instance \cite{openai2026math}). While the speed and carelessness of this delivery has been disorienting, these events motivated me to take a  backward look at one of my rewarding research directions. In this article I document an account of this journey. My goal in doing this is three fold: to provide the reader with some (very much non AI-generatable) context to my perspectives in this program; to communicate some simple thoughts on a near-future direction; and most importantly, to re-emphasize the importance of the human experience in mathematics. It should thence not come as a surprise that I did not use any AI to write this document. 

\subsection*{Leading to November 17th, 2022}

In the first few paragraphs, I explain how my interests moved towards the topic of selflessness. Between 2014 and 2018, I was an undergraduate student at UC Berkeley, where I wrote an honors thesis under the supervision of Dan Voiculescu and Brent Nelson (who was an NSF postdoctoral fellow at the time). My thesis project was to understand and provide an exposition of the seminal work \cite{VoiculescuFreeEntropy2}, where the framework for microstates free entropy dimension was established. My motivation for choosing this topic came from a talk given by Ben Hayes in the Probabilistic Operator Algebras Seminar in 2017 on his then preprint \cite{Hayes2018}. Prior to his talk, Brent Nelson explained to me a high level picture of free probability (notably free entropy \cite{DVVfreeentropy}) the notion of ``concentration of measure'', a'la Gromov-Milman's theorem. This is what got my attention, and since then I have harbored a fundamental interest in free entropy and free independence. As part of the project, I was directed by Voiculescu to read various papers of Kenley Jung, beginning with his key shift to the packing number perspective \cite{JungLemma}, and his work on tubularity \cite{JungTubularity}. The second of the two in particular fueled my interest and intuition in the concept of \emph{ultrapowers} of II$_1$ factors, especially in relation to free probability and entropy theory of II$_1$ factors. 

In the Fall of 2018 I joined the PhD program in Vanderbilt University to work with Jesse Peterson. In early 2019, Scott Atkinson (who was a postdoc at the time,) and I completed a paper \cite{ScottSri2019ultraproduct} that expanded some of the themes from \cite{JungTubularity}, notably studying amenability and conjugacy in the context of ultrapowers of II$_1$ factors. During the process of writing this article, I became interested in the concept of elementary equivalence, i.e, isomorphism of ultrapowers in the setting of II$_1$ factors. Upon discussing with Jesse Peterson, I learned that he too shared an interest in this concept, in particular the then open problem of whether there existed a pair of separable non Gamma II$_1$ factors that are not elementarily equivalent. This problem greatly captivated my interest too. He mentioned to me that he had worked on this problem with Andreas Thom earlier, in connection to Gabor Kun's work on sofic approximations of Kazhdan groups \cite{kun2016soficapproximationspropertytgroups}, which lead me to read more of Thom's papers about sofic groups and geometric and combinatorial group theory at large. Fortunately, I was able to take courses with Mark Sapir and Denis Osin on such topics and also interact greatly with Simon Andre (who was a postdoc at the time). In particular Simon explained to me Tarski's problem on elementary equivalence of free groups, and the ideas going into Sela's works \cite{sela1}. While I knew that the situation for the free group factors was rather hopeless at the time, I had developed a curiosity in possibly seeing if the solution to the Tarski's problem can extend to the setting of the $C^*$-ultrapower of the reduced free group $C^*$-algebras (i.e, $C^*$-elementary equivalence).

In the summer of 2022 as I was wrapping up my studies at Vanderbilt and moving to UCLA as a postdoc, I was fortunately able to make progress on the earlier thread of constructing a pair of non elementarily equivalent non Gamma II$_1$ factors (leading to the work \cite{exoticCIKE}). However, I was still interested in, and rather clueless about the $C^*$-algebra situation. Up until this point, I knew almost nothing about the program of classifying nuclear $C^*$-algebras (let alone non nuclear, and let alone up to elementary equivalence). Then came November 17th 2022, when Chris Schafhauser visited UCLA to give a department colloquium regarding quasidiagonal groups. I brought up this problem to Chris, and we had multiple discussions on an approach to solve it, albeit negatively (going in the opposite spirit of the group case). In this conversation, Chris educated me about the notion of strict comparison \cite{Blackadar_1989} of positive elements. To me this was a very natural property to arrive at in this context. Indeed, the way to prove conjugacy for hyperfinite subalgebras of ultraproducts of II$_1$ factors was using Murray von Neumann comparison, and strict comparison is a very reasonable replacement. In practice, strict comparison concretely allowed us to invoke some of Chris's intuition to access the $K_1$-group via computing exponential rank, as eventually demonstrated in \cite{elayavalli2025negativeresolutioncalgebraictarski}. 

Unfortunately, Chris told me that it was a dead end because strict comparison for $C^*_r(\mathbb{F}_n)$ was a well known open problem, and that the approach for $C^*_r(\mathbb{F}_\infty)$ from \cite{robert2012classification} had a clear obstruction to generalize. The obstruction was simply due to the fact that $C^*_r(\mathbb{F}_{\infty})$ had the nice property that any finite set can be perturbed a bit to locate an ``unused letter'' that is freely independent from the entire perturbed subalgebra. The main argument of Dykema-Rordam in \cite{dykema1998projections} was using the key free probabilistic fact: if $p,q$ are freely independent projections in a tracial $C^*$-algebra then $\|p-pqp\|< 1$. This then allows for building nicely located invertible elements and carrying out the necessary algebraic comparison.   Given the proximity to free probability, it was always on the back burner for me to try to understand strict comparison and potentially prove it for $C^*_r(\mathbb{F}_2)$. 

\subsection*{Towards selflessness of $C^*_r(\mathbb{F}_2)$} 

In the summer of 2023, I had the opportunity to travel to Europe. In particular I visited Munster and gave talks about proper proximality \cite{ding2023properly}. During this visit, I had an idea to attack the problem of strict comparison of $C^*_r(\mathbb{F}_2)$ based on a dynamical approach. It was to exploit the well known isomorphism on the group level itself $\mathbb{F}_2= \mathbb{F}_\infty\rtimes \mathbb{Z}$ (with the shift action). Unfortunately I could not make any progress despite this. I then lost sight of this problem until August 2024, when I was getting ready to prepare my applications for a tenure track position. I wrote to NC Phillips letting him know that I had applied to Oregon for a position, and he responded with a message to look at the strict comparison problem for $C^*_r(\mathbb{F}_2)$ since I had mentioned in my application materials that I worked on free probability. This piqued my interest, especially because I had invited Tattwamasi Amrutam a fellow Indian postdoc who worked on $C^*$-matters with Mehrdad Kalantar to visit UCSD for a month. 

The day he arrived, together with David Gao and Gregory Patchell (we were some sort of a trinity in UCSD, spent all our time doing math together) we sat down and decided to work on this problem. What is rather funny is that on day 2, we already solved the problem. What made all the difference was finding Leonel Robert's preprint on the notion of \emph{selflessness} from mid-2023 \cite{robertselfless}. Roughly speaking, this notion asks for a reduced freely independent copy of a Haar unitary in the $C^*$-ultrapower. To me, this rang many bells. First of all, in the II$_1$ factor context, this is extremely reminiscent one of my favorite theorems of Popa \cite{popa1995free}. Secondly, the group theoretic version of this is identical to the \emph{mixed identity free} property for groups. I had known about this property already from my times discussing with Simon Andre in Vanderbilt, in relation with one of my favorite lemmas of Baumslag which realized the surface groups as limit groups (which I put to use in the note \cite{sri}). It felt to me that a lot of things were coming together, and so I was very hopeful that we would solve this problem by using the mixed identity freeness of the free groups and somehow bumping it up to the $C^*$-algebra. After a failed attempt at taking Popa's proof and trying to adapt it, we came up with two approaches, one was to use Fell's absorption which we couldn't make work, and the other was to use random matrix theory. I was motivated to think in terms of random matrices because of developments in strong convergence, and my interest at that time in proving Peterson-Thom's conjecture for surface groups extending Ben's strategy \cite{HayesPT}. Fortunately, at the end of the day we could bootstrap a fine tuned version of mixed identity freeness using rapid decay, as demonstrated in the work \cite{loudermagee2025limitgroups}. The rest of the effort for our paper \cite{amrutam2025strictcomparisonreducedgroup} was to develop some expertise in geometric group theory (involving reading \cite{YANG_2014} and \cite{dahmani2017hyperbolically}) to push the idea as far as we could before publishing. Meanwhile I did not lose sight of my original purpose of all of this, so I emailed Chris rather excitedly and we had a meeting where I explained to him the proof. In parallel we resurrected our original plan and finished the paper settling the $C^*$-algebraic Tarski problem \cite{elayavalli2025negativeresolutioncalgebraictarski}. 

Soon after, Leonel Robert and I were interested in extending the approach of quantitative mixed identity freeness plus rapid decay, to see if we could get selflessness for families of reduced free products outside the group setting. Fortunately, Ben Hayes was already developing rapid decay for filtrations for other reasons involving the Peterson-Thom conjecture, so we were able to work together to produce a successful result for the free semicircular algebras \cite{HKER}. The number of selflessness results seemed to skyrocket around this time. Itamar Vigdorovich, who was in day to day mathematical contact with us at UCSD, was able to succesfully employ this strategy (thanks to Lafforgue's rapid decay result) in the case of linear groups leading to examples of selflessness among higher rank groups \cite{vigdorovich2025structuralpropertiesreducedcalgebras}. On the side of twisted reduced $C^*$-algebras, the work \cite{raum2025strictcomparisontwistedgroup} handled strict comparison despite having non-trivial center in the group which would obstruct selflessness. Together with Greg Patchell and a young Lizzy Teryoshin, we were able to push the free product results slightly beyond property RD, using Ricard-Xu's Khintchine inequality \cite{elayavalli2025remarks}. Thiel and co-authors were also obtaining various consequences of our result for the free groups to the finer structure of the $C^*$-algebra (for instance \cite{perera2025extensionspurecalgebras}). Notably, Greg Patchell, Leonel Robert, Itamar Vigdorovich and myself were involved in a very hard effort to capture $\mathbb{F}_2\times \mathbb{F}_2$ which clearly violated any MIF phenomena. We had developed a pretty nice harmonic analytic approach, but it was very technical and seemed hard to execute. Meanwhile, Ben Hayes, Greg Patchell, Leonel Robert and myself introduced the notion of selfless inclusions of $C^*$-algebras \cite{hayes2025selfless}, paralleling themes of $C^*$-irreducible inclusions or $\mathcal{Z}$-stable inclusions. It was quite applicable since we were able to exploit the position of the free unitary quantum group $C^*$-algebras between a selfless inclusion to conclude selflessness for examples that are not free products. Around this time I moved to Maryland and found myself rather stuck with no progress on selflessness.

\subsection*{Fallout from Toeplitz--Pimsner methods}
During a conference in Germany in late summer 2025, after listening to a talk of Greg Patchell about our work, Ozawa was inspired to find a new approach to prove selflessness \cite{ozawa2025proximalityselflessnessgroupcalgebras}. This exploited Pimsner's isomorphism theorem of universal Toeplitz algebras to reduced free products \cite{pimsner1997class}. In effect, the problem reduced to locating a Toeplitz operator in the $C^*$-ultrapower and check the Toeplitz relations (namely, $T^*aT=\rho(a)$  for every $a\in A$, where $\rho$ is the state on $A$) with the diagonal embedding. Ozawa made a designer combinatorial property called PHP for countable groups that allows to product this isometry. In fact he obtained a stronger property called \emph{complete selflessness} from this (that also bypasses the need for exactness in some arguments). This was used to vastly expand on the selflessness phenomena in countable groups. This includes arbitrary groups admitting extreme boundaries including arbitrary acylindrically hyperbolic groups \cite{yang2025extreme}, direct products with an exactness assumption, linear groups \cite{vigdorovich2026selflessreducedcalgebraslinear}, etc. This was a great boost to the program, enabling me to enter again with renewed interest. The first and rather elementary step was to move Ozawa's approach out of the group setting, which resulted in accessing selflessness for more general reduced free products \cite{fmmm2025selflessfreeprod, hayes2025selfless} and amalgamated free product $C^*$-algebras and HNN extensions \cite{gao2026selflessreducedamalgamatedfree}. We also went back and revised the selfless inclusions paper into the framework of the PHP property.  At this point, there were also various other applications of strict comparison for free free products, notably including an isomorphism theorem \cite{hirshberg2026isomorphism}.

Given all this progress, in late Fall 2025 my mind was set on two particular directions for moving the theory forward. The first direction was to revisit the early random matrix theory connections, especially the phenomena of strong convergence. In particular I was hoping that one could prove strong convergence of unitary representations or at least pieces of the argument by appealing to the new techniques in selflessness. The second direction was to expand the notion of selflessness beyond $C^*$-probability spaces, and to find the correct notion of relative selflessness.  To be more precise, let us say we have an inclusion of $C^*$-probability spaces $A\subset (B,\rho)$. If we wish to define a notion of ``$B$ is selfless relative to $A$'', we need to find replacements for two notions for the unitary $u\in B^\mathcal{U}$: the freeness, and the Haarness. Naively, one could suggest that $u$ is free with amalgamation from $A$, i.e, the $C^*$-algebra generated by $B$ and $u$ is isomorphic to the reduced amalgamated free product $B*_AC^*(A, u)$. As far as the second condition goes, one could suggest that $E_A(u^n)=0$ for every $n\in \mathbb{N}$. These two conditions make for a decent definition but the vagueness of the structure of the $C^*$-algebra generated by $A$ and $u$ meant that there was very little traction in developing this theory in a satisfying manner.

To address these two matters (and meet my friends again), I decided to organize a conference at the University of Maryland (thanks to the generosity of the Brin Mathematical Research Center) in January 2026. This turned out to be one of the most enjoyable and invigorating research experiences, wherein new directions were pushed. First of all, together with co-authors David Gao, Marius Junge, Greg Patchell and Leonel Robert, we managed to identify the missing piece to develop a satisfying theory of relative selflessness: enforcing that $u$ commutes with $A$! This was extremely natural for multiple reasons, notably including the compatibility with the structure of the isomorphism arising from Pimsner's theorem for conditional expectations. This allowed us to push even further and develop a satisfying theory of selflessness for arbitrary $C^*$-correspondences \cite{gao2026selfless} with various applications. On a parallel side, David Gao and I found a conceptually new proof of Pimsner-Voiculescu's non isomorphism theorem \cite{gao2026nonisomorphism}, avoiding $K$-theory and using a brand new idea of considering embedding spaces in II$_1$ factors and locating many freely independent haar unitaries. Then approximately a month later I visited Austin for a conference on sofic groups organized by Lewis Bowen. In this workshop, David Gao, myself, Aareyan Manzoor and Greg Patchell managed to find a desirable connection back to random matrix theory \cite{gao2026newsourcepurelyfinite}. By developing a new inductive approach weaving together $C^*$-ultrapowers and combinatorial group theory we were able to exploit Pimsner's isomorphism theorem and prove strong convergence of unitary representations, unifying and significantly generalizing prior works on the MF and PMF property. Our approach was also able to access strong convergence for unitary representations with finite images of 3-manifold groups, a property of interest to geometers \cite{song2025randomharmonicmapsspheres}. A particularly useful tool motivated by this result that David Gao and I emphasized is Toeplitz exactness \cite{gao2026toeplitzexactnessstrongconvergence}. The last result I wish to comment on is Ozawa's powerful uniform exactness theorem \cite{ozawa2026uniform} for free groups, which appeared to be motivated by these prior strong convergence results and selflessness. Not only does this open up a new surprising phenomena, but also offers a new proof of selflessness for $C^*_r(\mathbb{F}_2)$, without needing rapid decay.

\subsection*{Looking forward}
Since the beginning of mid August 2026, there has been a further outburst of results on selflessness via AI use. I will not outline each and every one of them but emphasize one crucial result \cite{robert2026strictcomparisonselflessness} that provided closure to the question of the relationship between selflessness and strict comparison. Given the fact that Popa's theorem \cite{popa1995free} is clearly a characterization of II$_1$ factors, there was a suspicion among experts that these two ought to be equivalent, under additional necessary assumptions including simplicity and unique trace. This was indeed proven to be the case, and the proof interestingly uses Popa's theorem and the strict comparison hypothesis to adjust the Popa unitary to put it in position to apply Pimsner's isomorphism theorem again to verify the $C^*$-freeness. One very important consequence of this idea is a new proof of the Blackadar-Kirchberg MF property \cite{Blackadar1997Generalized} for $C^*_r(\mathbb{F}_2)$, which was originally proved by Haagerup and Thorbjornsen \cite{HaagerupThorbjornsen2005}. This interested me quite a lot, especially since Lizzy Teryoshin and myself were working towards this when the announcement came out. It is also important to point out that one of the main open questions in this theory still remains open: is every $C^*$-simple group $C^*$-selfless. It would certainly be interesting to understand the passage from simplicity to selflessness in the group case (without being able to access strict comparison apriori). 

 By now, the parallel between tracial selfless $C^*$-algebras and the class of II$_1$ factors among von Neumann algebras, is very clear. In the case of arbitrary II$_1$ factors, Popa has uncovered quite a lot of internal structure beyond just his free independence result. For instance in Popa's theorem already, it is shown that if $N\subset M$ is an irreducible subfactor, then the inclusion is selfless in the II$_1$ factor sense, i.e, one can locate a freely independent haar unitary in the II$_1$ factor ultrapower of $N$. In this context, an important early theorem of Popa \cite{Popa1981} with various applications  is the existence of an irreducible copy of $R$ inside any II$_1$ factor. Fusing these two philosophies, I am tempted to make the following natural conjecture. Note that Robert has already proved the basic statement that selfless $C^*$-algebras admit a copy of $\mathcal{Z}$ \cite{robert2012classification}.

\begin{conjecture}
    Every selfless $C^*$-algebra admits a selfless inclusion (in the sense of \cite{hayes2025selfless}) of the Jiang-Su algebra $\mathcal{Z}$.
\end{conjecture}

Pushing forward a couple of decades, Popa proved a very useful strengthening of his earlier irreducible embedding theorem: every II$_1$ factor admits a ``coarse'' hyperfinite subfactor \cite{Popa2021Coarse}. This is a much stronger statement, and has already been used nicely (for instance to construct very exotic II$_1$ factors that do not arise as a non trivial crossed product \cite{quero2026classii1factorsnontrivial}). Therefore it is very reasonable to expect the same phenomena in selfless $C^*$-algebras, but I do not know what is the ``correct'' analogy of coarseness here. In this same vein, another theorem that has been proved globally for all II$_1$ factors again using Popa's theorem is single generation. While this seems to have been already prompted in \cite{elliott2026generatorsseparableunitalcalgebras} I am definitely interested in understanding single generation in selfless $C^*$-algebras at some point. Yet another result about all II$_1$ factors using Popa's theorem is Jekel's contractability theorem \cite{Jekel2025SOT}. I wonder if there is a ``correct'' analogy for selfless $C^*$-algebras here too. On a grander level, I think the emergence of selflessness might allow for more import of progress in the rigidity and classification of nonamenable II$_1$ factors, into the regime of non-nuclear $C^*$-algebras.

Regarding existing developments concerning selflessness, there are some  loose ends that appear to remain untied, which I point out here. Popa's relative Dixmier property paper \cite{Popa2000_Dixmier} doesn't appear to have a selfless analogue yet. Namely, is it true that if $H<G$ are two icc groups, then $G$ is $C^*$-selfless iff $H$ is $C^*$-selfless? Another problem is concerning tensor products.
    Suppose $(A_i,\rho_i)$, $i=1,2$ are  $C^*$-probability spaces, and $(A_1,\rho _1)$ is selfless, and $(A_2,\rho_2)$ is simple with unique trace. Then is the minimal tensor product $(A_1\otimes A_2, \rho_1\otimes \rho_2)$ selfless?  The above is known if one additionally assumes $A_1$ is exact \cite{gao2026selfless}. In this spirit, the more general situation for crossed products also appears open:    Let $G$ be a countable $C^*$-selfless group that acts on a simple monotracial $C^*$-algebra $A$. Then is the reduced crossed product $A\rtimes_r G$ selfless?    The above is known if $G$ has the PHP property, and the action is approximately inner \cite{gao2026selfless}. 
    Additionally, concerning the PHP property, for a countable group $G$, it also seems unknown if PHP is equivalent to complete selflessness of $C^*_r(G)$. Regarding free products, the following is a natural question: among the class of reduced free products of non trivial $C^*$-probability spaces, is simplicity equivalent to selflessness? 
In the case that the $C^*$-probability spaces are finite dimensional, the above has been settled in the upcoming work of Gould, Hartglass and Robert {\cite{GouldHartglassRobert}}. More generally,
    suppose $(A_i,\rho_i)$ are $C^*$-probability spaces, and a $B\subset A_i$ with Watatani index greater than 2 for each $i=1,2$. Then is the amalgamated free product $A_1*_BA_2$ selfless if and only if it is simple? What about a classification of strict comparison? It seems that strict comparison should be more easily accessible outside the simple and unique trace examples.

    \subsection*{A concluding remark}
    With the reader's permission, I would like to share that personally I am of the following mindset: despite the recent emergence of AI tools that vastly accelerate the production of theorems and proofs, it would be beneficial if our community were to perhaps slow down and re-establish a grip on the mathematics we collectively cherish. I believe we should not opt for ``instant gratification'' (prompt and claim credit for AI assisted results in quick succession), but instead focus on ``delayed gratification'' (fostering genuine human mathematical connection and understanding).

\subsection*{Acknowledgements} 
I owe a special acknowledgment to Leonel Robert. His beautiful contributions to the theory of $C^*$-algebras have inspired and enabled rewarding mathematical exploits of myself and many other mathematicians. At every turn, starting from our first correspondence, Leonel showed me nothing but generosity, patience, kindness and camaraderie. He has been a  source of inspiration as a mathematician and a human being and it is my fortune to call him friend. I am grateful to Tattwamasi Amrutam, David Gao, Ben Hayes, Marius Junge, Aareyan Manzoor, Mahan Mj, Greg Patchell, Leonel Robert, Chris Schafhauser, Lizzy Teryoshin, Itamar Vigdorovich for sharing projects with me in this $C^*$-algebraic journey in the past few year years. I am also grateful to the following kind people who have supported me greatly through a painful recovery from a near fatal paragliding crash in September 2026: James Chesley, Ivo Duijnstee, Karthik Ganapathy, David Gao, Mike Hartglass, Mary Reed Johnson Cindy Looy, Nicki Magill, Andrew Marks, Tejas Nareddy, Krishna Parthasarathy, Varun Pattabiraman, Michael Ragone, Bharatram Rangarajan, Marc Rieffel, Brandon Seward, Forte Shinko, Lizzy Teryoshin and Jennifer Zhu.  Last but not the least, I am grateful to the  following puppies who have blessed me with true love along my journey: Scampy Bejleri, Martok Castillo, Jiu-Jiu Frank, Evie Hartglass, Elsa Henderson, Olive Johnson, Jeeva Muthukumar, Kramer Nelson, Bandit Parthasarathy, Murphy Patchell, Wolfsbane Peterson, Fuffle Sigurdsen, Kona Soltani, Rocky Thiruvengadam, Runa Thom, Manju Varadhan, Rover Vijayaraghavan, Bailey Weiwu, and Duchess Whiteaker.

\bibliographystyle{alpha}
	\bibliography{references}

\end{document}